\documentclass[11pt]{amsart}

\usepackage{amscd}
\usepackage{amsfonts}
\usepackage{amssymb}
\usepackage{amsmath}
\numberwithin{equation}{section}

\newtheorem{pro}{Proposition}[section]
\newtheorem*{t*}{{\bf Theorem}}

\newtheorem{definition}[pro]{Definition}
\newtheorem{remark}[pro]{Remark}

\newcommand{\be}{\beta}
\newcommand{\ga}{\gamma}

\newcommand{\om}{\omega}

\newcommand{\ra}{{\rightarrow}}
\newcommand{\lra}{{\longrightarrow}}

\newcommand{\ie}{{i.e$.$\,}}
\newcommand{\cf}{{cf$.$\,}}

\newcommand{\CC}{{\mathbb C}}

\newcommand{\HH}{{\mathbb H}}

\newcommand{\PP}{{\mathbb P}}
\newcommand{\RR}{{\mathbb R}}

\newcommand{\ZZ}{{\mathbb Z}}

\newcommand{\Iso}{\mathop{\rm Iso}\nolimits}

\newcommand{\Null}{\mathop{\rm Null}\nolimits}

\newcommand{\Aut}{\mathop{\rm Aut}\nolimits}

\begin{document}
\pagestyle{myheadings} \thispagestyle{empty} \setcounter{page}{1}

\title[\ ]
{Nonexistence of cusp cross-section of one-cusped complete complex
hyperbolic manifolds II}
\author[\ ]{Yoshinobu KAMISHIMA}
\address{Department of Mathematics, Tokyo Metropolitan
University, Minami-Ohsawa 1-1, Hachioji, Tokyo 192-0397, JAPAN}
\email{kami@comp.metro-u.ac.jp} 
\date{\today}

\markboth{Yoshinobu KAMISHIMA}{Cusp representation}

\begin{abstract}
Long and Reid have shown that some compact flat 3-manifold cannot be
diffeomorphic to a cusp cross-section of any complete finite volume
$1$-cusped hyperbolic $4$-manifold. Similar to the flat case, we
give a negative answer that there exists a $3$-dimensional closed
Heisenberg infranilmanifold whose diffeomorphism class cannot be
arisen as a cusp cross-section of any complete finite volume
$1$-cusped complex hyperbolic $2$-manifold. This is obtained from
the formula by the characteristic numbers of bounded domains related
to the Burns-Epstein invariant on strictly pseudo-convex
$CR$-manifolds  \cite {b-he},\cite{b-e2}. This paper is a sequel of
our paper\cite{ka3}. \end{abstract}

\maketitle


\setcounter{section}{1}

\section*{Introduction}\label{INT}
We shall consider whether every Heisenberg infranilmanifold can be
arisen, up to diffeomorphism, as a cusp cross-section of a complete
finite volume $1$- cusped complex hyperbolic manifold. Long and Reid
considered the problem that
 every compact Riemannian flat manifold is diffeomorphic to
a cusp cross-section of a complete finite volume $1$-cusped
hyperbolic manifold. They have shown it is false for some compact
flat $3$-manifold \cite{l-r1}. We shall give a negative answer
similarly to the flat case.

\begin{t*}\label{th1}
Any $3$-dimensional closed Heisenberg infranilmanifold with
nontrivial holonomy cannot be
diffeomorphic to a cusp cross-section of any complete finite volume
$1$-cusped complex hyperbolic $2$-manifold.
\end{t*}
\noindent  McReynolds informed us that W. Neumann and A. Reid have
obtained the similar result.

\section{Heisenberg infranilmaniold}
Let  $\displaystyle \langle z,w \rangle=\bar z_1\cdot w_1+\bar
z_2\cdot w_2 +\cdots+\bar z_n\cdot w_n$ be the Hermitian inner
product defined on $\CC^n$. The Heisenberg nilpotent Lie group
$\mathcal N$ is the product $\RR\times \CC^n$ with group law:
\begin{equation}\label{product}
(a,z)\cdot (b,w)=(a+b-{\rm Im}\langle z,w\rangle,\ z+w).
\end{equation}
It is easy to see that $\mathcal N$ is $2$-step nilpotent, \ie
$[\mathcal N,\mathcal N]=(\RR,0)=\RR$, which is the central subgroup
$\mathcal C(\mathcal N)$ of $\mathcal N$. This induces a central
group extension: $\displaystyle 1\ra \mathcal C(\mathcal N)\ra
\mathcal N \stackrel{P}\lra \CC^n\ra 1$. Let
$\mathop{\Iso}(\HH\sb{\CC}\sp {n+1})$ be the full group of the
isometries of the complex hyperbolic space $\HH\sb{\CC}\sp {n+1}$.
It is isomorphic to ${\rm PU}(n+1,1)\rtimes \langle\tau \rangle$
where $\tau$ is the (anti-holomorphic) involution induced by the
complex conjugation.
 The Heisenberg rigid motions is defined as a subgroup of the
stabilizer $\mathop{\Iso}(\HH\sb{\CC}\sp {n+1})_\infty$ at the point
at infinity $\infty$ .

\begin{definition}\label{gene-nil}The group of Heisenberg rigid motions
$\displaystyle{\rm E}^{\tau}(\mathcal N)$ is defined to be $\mathcal
N\rtimes ({\rm U}(n)\rtimes \langle \tau\rangle)$. A Heisenberg
infranilmanifold (respectively orbifold) is a compact manifold
(respectively orbifold) $\mathcal N/\pi$ such that $\pi$ is a
torsionfree (not necessarily torsionfree) discrete cocompact
subgroup of ${\rm E}^{\tau}(\mathcal N)$.
\end{definition}

\section{$CR$-structure on $S^{2n-1}-S^{2n-1}$}\label{sphere-contact}
 The sphere complement $S^{2n+1}-S^{2n-1}$ is
a spherical $CR$ manifold with the transitive group
$\mathop{\Aut}_{CR}(S^{2n+1}-S^{2n-1})$ of $CR$ transformations
which is isomorphic to the unitary Lorentz group ${\rm U}(n,1)$.
Note that $S^{2n+1}-S^{2n-1}$ is identified with the
$(2n+1)$-dimensional Lorentz standard space form $V_{-1}^{2n+1}$ of
constant sectional curvature $-1$. The center $\mathcal Z{\rm
U}(n,1)$ of ${\rm U}(n,1)$ is $S^1$. Then $V_{-1}^{2m+1}$ is the
total space of the principal $S^1$-bundle over the complex
hyperbolic space: $\displaystyle S^1\ra ({\rm
U}(n,1),V_{-1}^{2n+1})\stackrel{\nu}\lra ({\rm
PU}(n,1),\HH_{\CC}^n)$. If $\omega_{\HH}$ is the connection form of
the above principal bundle, then it is a contact form on
$V_{-1}^{2n+1}$ such that $\mathop{\Null}\om_{\HH}$ is a $CR$
structure. Note that $d\om_\HH=\nu^*\Omega_\HH$ up to constant
factor for the K\"ahler form $\Omega_\HH$ on $\HH_{\CC}^n$. Since
${\rm U}(n,1)=S^1\cdot {\rm SU}(n,1)$, the above equivariant
principal bundle induces the following commutative fibrations:
\begin{equation}\label{SLorentz}
\begin{CD}
\ZZ@>>> (\hat{{\rm SU}}(n,1), \tilde V_{-1}^{2n+1})@>\hat\nu>> ({\rm
PU}(n,1),\HH_{\CC}^n) \\
@VVV     @VVV              {||} @.\\
\ZZ/n+1 @>>> ({\rm SU}(n,1), V_{-1}^{2n+1})@> \nu >> ({\rm
PU}(n,1),\HH_{\CC}^n).
\end{CD}
\end{equation}Here $\displaystyle \hat{{\rm SU}}(n,1)$ is a lift of
${\rm SU}(n,1)$ associated to the covering $\ZZ\ra \tilde
V_{-1}^{2n+1}\ra  V_{-1}^{2n+1}$.
 For a discrete
subgroup $G\subset {\rm PU}(n+1,1)$ such that $\HH^{n+1}_{\CC}/G$ is
a complete finite volume complex hyperbolic orbifold, let $\hat
G\subset{\hat{\rm SU}}(n,1)$ be a lift where $1\ra \ZZ\ra \hat G\ra
G\ra 1$ is an exact sequence. Then $\displaystyle S^1\ra \tilde
V_{-1}^{2n+1}/\hat G\stackrel{\hat\nu} \lra \HH^{n+1}_{\CC}/G$ is an
injective Seifert fibration (\ie the singular fiber bundle with
typical fiber is $S^1$. The exceptional fiber is also a circle.)

\section{Burns and Epstein's formula}\label{Bu-Ep} In general, the Heisenberg
infranilmanifold or its two fold cover at least admits a spherical
$CR$-structure, see Definition \ref{gene-nil}. In \cite{b-e1}, Burns
and Epstein obtained the $CR$-invariant $\mu(M)$ on the
$3$-dimensional strictly pseudoconvex $CR$-manifolds $M$ provided
that the holomorphic line bundle is trivial. Let $X$ be a compact
strictly pseudoconvex complex $2$-dimensional manifold with smooth
boundary $M$. Then they have shown the following equality in
\cite{b-e2}:
\begin{equation}\label{eq3}
\int_X c_2-\frac 13 c_1^2=\chi(X)-\frac 13 \int_X\bar c_1^2+\mu(M).
\end{equation}
Here $\bar c_1$ is a lift of $c_1$ by the inclusion
$j^*:H^2(X,M;\RR)\ra H^2(X;\RR)$.

\section{Geometric boundary}
\subsection{One-cusped complex hyperbolic $2$-manifold}
Let ${\rm E}^\tau(\mathcal N)$ be the group of Heisenberg rigid
motions on the $3$-dimensional Heisenberg nilpotent Lie group
$\mathcal N$ and $L:{\rm E}^\tau(\mathcal N)\ra {\rm U}(1)\rtimes
\langle\tau\rangle$ the holonomy homomorphism. Suppose that
$M=\mathcal N/\Gamma$ is realized as a cusp cross-section of a
complete finite volume one-cusped complex hyperbolic $2$-manifold
$W=\HH^2_\CC/G$. Put $\bar W=\HH^2_\CC/G-M\times (0,\infty)$ so that
$\partial {\bar W}=M$. Then $\bar W$ is homotopic to $W$ and $M$ is
viewed as a boundary of ${\rm Int}\bar W$ which supports a complete
complex hyperbolic structure. The holonomy group $L(\Gamma)$ of a
$3$-dimensional compact Heisenberg non-homogeneous infranilmanifold
$M=\mathcal N/\Gamma$ is a cyclic subgroup of order $2, 3,4, 6$ of
${\rm U}(1)$ or $L(\Gamma)$ is isomorphic to $\ZZ/2\times
\ZZ/2\subset{\rm U}(1)\rtimes \langle\tau\rangle$, see \cite{dekim},
\cite{ben} for the classification. If $M$ has the holonomy
$\ZZ/2\times \ZZ/2$, then $G$ has nontrivial summand in
$\langle\tau\rangle$ of ${\rm Iso}(\HH^2_\CC)={\rm PU}(2,1)\rtimes
\langle\tau\rangle$. The two fold cover $W/G\cap{\rm PU}(2,1)$ is
still a one-cusped complex hyperbolic manifold for which the cusp
cross-section is the two fold cover of $M$ whose holonomy group
becomes $\ZZ/2\subset{\rm U}(1)$. When the holonomy group belongs to
${\rm U}(1)$, the spherical $CR$-structure on $M$ is canonically
induced from the complex hyperbolic structure on $W$. (Note that
$\tau$ does not preserve the $CR$-structure bundle.)

\subsection{Integral of ${\bf \bar c_1^2}$ }\label{integrand}
 Let $p:\tilde {W}\ra {W}$ be the finite
 covering, say of order $\ell$, whose induced covering $\tilde M$ of $M$
is now a (homogeneous) nilmanifold (using the separability argument
if necessary). Possibly it consists of a finite number of such
nilmanifolds. Since $W$ admits a complete Einstein-K\"ahler metric,
we know that $\displaystyle c_2-\frac 13 c_1^2=0$. Moreover, since
$\tilde M$ is a spherical $CR$ manifold with trivial holomorphic
line bundle, it follows that $\mu(\tilde M)=0$.
 As in $\S \ref{Bu-Ep}$, let $j^*:H^2(\bar W,M:\RR)
\ra H^2(\bar W:\RR)=H^2(W:\RR)$
 be the map such that $j^*\bar c_1(\bar W)=c_1(W)$.
Applying \eqref{eq3} to $\tilde {W}$, we have $\displaystyle
\chi(\tilde W)=\frac 13 \int_{\tilde W}\bar c_1^2$. As $p^*(\bar
c_1(\bar W))=\bar c_1(\tilde W)$ by naturality and $p_*[\tilde W]
=\ell[\bar W]$, it follows that $\displaystyle \int_{\tilde W}\bar
c_1^2=\langle \bar c_1^2(\tilde W),[\tilde W]\rangle=\langle \bar
c_1^2(\bar W),\ell[\bar W]\rangle$. Since $\chi(\tilde
W)=\ell\chi(W)$, $\displaystyle 3\chi( W)=\langle \bar c_1^2(\bar
W),[\bar W]\rangle$.
\begin{pro}\label{integervalue}
If $M=\mathcal N/\Gamma$ is realized as a cusp
cross-section of a complete finite volume one-cusped complex
hyperbolic $2$-manifold $W=\HH^2_\CC/G$, then
$\bar c_1^2(\bar W)$ is an integer in $H^4(\bar W,M:\ZZ)=\ZZ$.
\end{pro}

\subsection{Torsion element in $M$}\label{torsion of b}
Given a $CR$-structure on $M$, there is the canonical splitting
$TM\otimes \CC=B^{1,0}\oplus B^{0,1}$ where $B^{1,0}$ is the
holomorphic line bundle.  Since $M$ is an infranilmanifold but not
homogeneous, $B^{1,0}$ is nontrivial, \ie $c_1(B^{1,0})\neq 0$. (In
fact, it is a torsion element in $H^2(M:\ZZ)$,
 because the $\ell$-fold covering $\tilde M$ has
the trivial holomorphic bundle.) The spherical $CR$ manifold $M$ has
a characteristic $CR$ vector field $\xi$.  If $\epsilon^1$ is the
vector field on $M$ pointing outward to $W$, then the distribution
$\langle\epsilon^1,\xi \rangle$ generates a trivial holomorphic line
bundle $T\CC^{1,0}$ on $M$ for which $TW\otimes \CC|M=B^{1,0}+
T\CC^{1,0}\oplus B^{0,1}+T\CC^{0,1}$. As $\displaystyle
i^*(c_1(W))=c_1(B^{1,0}+T\CC^{1,0})=c_1(B^{1,0})$ and $\ell\cdot
c_1(B^{1,0})=0$,
 we have $j^*\be=\ell\cdot c_1(W)$ for some integral class $\be\in H^2(\bar
W,M:\ZZ)$.
\subsection{${\bf H_1(M:\ZZ)}$}\label{H1}
 Let $\displaystyle 1\ra\Delta\ra
\Gamma\ra F\ra 1$ be the group extension of the fundamental group
$\Gamma=\pi_1(M)$ where $\Delta$ is the maximal normal nilpotent
subgroup and $F\cong \ZZ_\ell \ (\ell=2,3,4,6)$ or $F\cong
\ZZ_2\times \ZZ_2$. Recall that $\Delta$ is generated by $\{a,b,c\}$
where $[a,b]=aba^{-1}b^{-1}=c^k$ for some $k>0$. It follows that
$\displaystyle \Delta/[\Delta,\Delta]= \ZZ\oplus\ZZ\oplus \ZZ_k$.
Let $\ga$ be an element of $\Gamma$ which maps to a generator of
$\ZZ_\ell$. A calculation shows that  {\rm mod} $[\Delta, \Delta]$,
\begin{equation}\label{cyclicrep}
\begin{split}
\ga a\ga^{-1}&=a^{-1}, \ \ga b\ga^{-1}=b^{-1}\ (\ell =2),\\
 \ga a\ga^{-1}&=b, \ \ \ \ \ga b\ga^{-1}=a^{-1}b^{-1}\ (\ell =3),\\
\ga a\ga^{-1}&=b, \ \ \ \ \ga b\ga^{-1}=a^{-1}\ (\ell=4),\\
 \ga a\ga^{-1}&=b, \ \ \ \ \ga b\ga^{-1}=a^{-1}b\ (\ell=6).
\end{split}\end{equation}
When $F=\ZZ_2\times\ZZ_2$, let $\delta$ be an element of $\Gamma$
which goes to another generator of $F$. Then $\ga a\ga^{-1}=a,\ \ga
b\ga^{-1}=b^{-1}$ {\rm mod} $[\Delta, \Delta]$. In view of the above
relation \eqref{cyclicrep}, $\ga$ (also $\delta$) becomes a torsion
element of order $m$ in $\Gamma/[\Gamma,\Gamma]$ where $m$ is
divisible by $\ell$. As $\Gamma$ is generated by $\{a,b,c,\ga\}$ or
$\{a,b,c,\ga,\delta\}$, it follows that
\begin{equation}\label{homorep}
\begin{split}
H_1(M:\ZZ)=\ZZ_k\oplus\ZZ_m\oplus\left\{\begin{array}{lr}
\ZZ_2\oplus\ZZ_2& (\ell=2)\\
\ZZ_3& (\ell=3)\\
\ZZ_2 & (\ell=4)\\
1& (\ell=6)\\
\ZZ_2\oplus\ZZ_2&
\end{array}\right\}.
\end{split}\end{equation}
In any case, if $\mathcal N/\Gamma$ has a nontrivial holonomy group
$F$, then $H_1(M:\ZZ)$ is a torsion group.

\subsection{Intersection number}\label{intersect}
 Put $\bar H^2(\bar W,M:\ZZ)=H^2(\bar W,M:\ZZ)/{\rm
 Tor}$ where ${\rm Tor}$ is the torsion subgroup.
We have a nondegenerate inner product $\displaystyle \bar H^2(\bar
W,M:\ZZ)\times \bar H^2(\bar W,M:\ZZ)\ra \ZZ$ defined by the
intersection form
\[ (x,y)=\langle x\cup y,[\bar W]\rangle.
\]Denote by $\bar W\#\pm\CC\PP^2$ the connected sum of
$\bar W$ with $\CC\PP^2\#-\CC\PP^2$ if necessary. We can assume that
$(\ , )$ is an indefinite form of odd type, \ie there are nonzero
elements $x, y\in \bar H^2(\bar W\#\pm\CC\PP^2,M:\ZZ)= \bar H^2(\bar
W,M)+\langle 1\rangle + \langle -1\rangle$ such that $(x,x)$ is odd
and $(y,y)=0$.  By $\langle \pm 1\rangle$ we shall mean that it is
generated by either $x_+$ or $x_{-}$ of $\bar H^2(\bar
W\#\pm\CC\PP^2,M:\ZZ)$ such that $(x_{\pm},x_{\pm})=\pm 1$
respectively. Moreover by the classification of nondegenerate
indefinite inner product \cf \cite{h-m}, there is an isomorphism
preserving the inner product from $\bar H^2(\bar
W\#\pm\CC\PP^2,M:\ZZ)$ onto
\begin{equation}\label{pre-iso}
 m\langle 1\rangle \oplus n\langle -1\rangle=
\langle 1\rangle_1\oplus\cdots\oplus \langle 1\rangle_m\oplus
\langle -1\rangle_{1}\oplus\cdots \oplus \langle -1\rangle_n
\end{equation}for $(m,n\neq 0)$. Here $\langle \pm 1\rangle_i$ is the $i$-th copy of
$\langle \pm 1\rangle$.
 Consider the commutative diagram:
\begin{equation}\label{poincare}
\begin{CD}
 H^2(\bar W,M:\ZZ)@>j^*>> H^2(\bar W:\ZZ)@>i^*>> H^2(M:\ZZ)\\
           @VDVV        @VDVV                 @VDVV \\
H_2(\bar W:\ZZ) @>j_*>> H_2(\bar W,M:\ZZ)@>\partial>> H_1(M:\ZZ).
\end{CD}
\end{equation}
It follows from \eqref{homorep} that $j_*:\bar H_2(\bar W:\ZZ)\ra
\bar H_2(\bar W,M:\ZZ)$ is injective and is isomorphic if $\ZZ$
replaces $\RR$. Similarly note that $\displaystyle j_*: \bar
H_2(\bar W\#\pm\CC\PP^2:\ZZ)\lra \bar H_2(\bar
W\#\pm\CC\PP^2,M:\ZZ)$ is injective (and an isomorphism for the
coefficient $\RR$). Identified the generators of $\bar H_2(\bar
W\#\pm\CC\PP^2:\ZZ)$ with the basis \eqref{pre-iso} of $\bar
H^2(\bar W\#\pm\CC\PP^2,M:\ZZ)$, we may choose the generators
$[V_i]\in \bar H_2(\bar W\#\pm\CC\PP^2,M:\ZZ)$ such that
\begin{equation}\label{gene}
j_*(\langle \pm 1\rangle_i)= \ell_i[V_i]
\end{equation}for some $\ell_i\in\ZZ$.

\subsection{Canonical bundle}\label{cano}
The circle (line) bundle $L$: $S^1\ra \tilde V_1/\hat G\ra
\HH^2_\CC/G=W$ is represented by the K\"ahler form $\Omega$ of the
K\"ahler-hyperbolic metric, \ie $[\Omega]=c_1(L)\in H^2(W:\ZZ)$.
Hence $W=\HH^2_\CC/G$ is projective-algebraic, \ie $W \subset
\CC\PP^N$ so $c_1(W)$ can be represented by $c_1([V])$ for some
divisor $V$ in $W$, \ie $D(c_1(W))=[V]\in \bar H_2(\bar W,M:\ZZ)$,
compare \cite{g-h}. Embed $V$ into $\bar W\#\pm\CC\PP^2$ and suppose
that
\[ [V]=\sum_{i}a_i[V_i]\in
\bar H_2(\bar W\#\pm\CC\PP^2,M:\ZZ).
\]
As $D\circ i^*c_1(W)=\partial [V]$, it follows $\ell\partial([V])=0$
by the argument of $\S$ \ref{torsion of b}.
 We observe that $\partial [V]$ maps into $\ZZ_m$ in $H_1(M:\ZZ)$
(\cf \eqref{homorep}) and so does each $\partial [V_i]$. It may
occur that $\partial a_i[V_i]=\partial a_j[V_j]$ for some $i,j$. So
we can write $[V]=k[V_1]+j_*x$ where $x\in H_2(\bar
W\#\pm\CC\PP^2:\ZZ)$ and $V_1$ satisfies that
\begin{itemize}
\item[(1)] $\partial {\bar V}_1=S^1$ and $\ell[S^1]=0$ in $\ZZ_m\subset H_1(M:\ZZ)$.
\item[(2)]  $\ell$ is minimal with respect to (1).
\item[(3)] $(k,l)$ is relatively prime.
\end{itemize}

\subsection{Realization of ${\bf \bar c_1}$}\label{first c-rep}
As $\ell\partial[V_1]=0$ in $H_1(M:\ZZ)$, there is a surface $U$ in
$W$ whose cycle $[U]\in H_2(\bar W\#\pm\CC\PP^2:\ZZ)$ represents
$j_*[U]=\ell[V_1]$.

Let $[U]=a_1\langle \pm 1\rangle_1 +a_2\langle \pm 1\rangle_2+
\cdots +a_s\langle \pm 1\rangle_s$. Then, $\ell
[V_1]=a_1\ell_1[V_1]+a_2\ell_2[V_2]+\cdots+a_s\ell_s[V_s]$. Since
each $[V_i]$ is a generator of $\bar H_2(\bar
W\#\pm\CC\PP^2,M:\ZZ)$, it follows that $\ell=a_1\ell_1$ and $a_j=0$
$(j\neq 1)$. Hence $[U]=a_1\langle \pm 1\rangle_1$. On the other
hand, note that $\langle \pm 1\rangle_1$ is a cycle of $\bar
H^2(\bar W\#\pm\CC\PP^2,M:\ZZ)$ for which $j_*(\langle \pm
1\rangle_1)=\ell_1[V_1]$ by \eqref{gene}. Noting that $\ell$ is
minimal by Property (2) of $\S$ \ref{cano}, $\ell_1$ is divisible by
$\ell$. Therefore $\ell_1=\pm\ell$ and $a_1=\pm 1$ so that
$[U]=\pm\langle \pm 1\rangle$. In particular, the intersection
number
\begin{equation}\label{odd1}
[U]\cdot [U]=\pm 1.
\end{equation}

Put $\displaystyle y=\frac k\ell[U]+x\in H_2(\bar
W\#\pm\CC\PP^2:\RR)$. Calculate
\begin{equation*}\begin{split}
y\cdot y&=\frac {k^2}{\ell^2}[U]\cdot [U] +\frac{2k}\ell[U]\cdot x+x\cdot x\\
        &= \pm\frac {k^2}{\ell^2} +\frac{2k}\ell[U]\cdot x+x\cdot x,
\end{split}
\end{equation*}

\begin{equation}\label{selfnum}
\ell (y\cdot y)=\pm\frac {k^2}{\ell}\ \ {\rm mod}\  \ZZ
\end{equation}Noting that $(k,l)=1$ by Property (3) of $\S$ \ref{cano},
if $\ell\neq 1$, $y\cdot y$ cannot be an integer.

As $\displaystyle j_*(\frac k\ell[U])=k[V_1]$, note that
$j_*y=k[V_1]+j_*x=[V]$. Consider the following diagram:
\begin{equation}\label{poincare1}
\begin{CD}
\bar H^2(\bar W\#\pm\CC\PP^2:\RR)@>j^*>> \bar H^2(\bar W\#\pm\CC\PP^2,M:\RR)\\
      @|      @|                \\
\bar H_2(\bar W:\RR)+ \langle 1\rangle + \langle -1\rangle @>j_*+
{\rm id}>> \bar H_2(\bar W,M:\RR)+\langle 1\rangle+ \langle
-1\rangle.
\end{CD}
\end{equation}Let $y=y_0+ t\langle 1\rangle+ s\langle
-1\rangle$ for some $y_0\in\bar H_2(\bar W:\RR), s,t\in \RR$. As
$j_*y=[V]$, it follows that $[V]=j_*y_0+t\langle 1\rangle+ s\langle
-1\rangle$. Noting $[V]\in \bar H_2(\bar W,M:\ZZ)$, we have that
$[V]=j_*y_0$ and $t=s=0$. In particular, this implies that
$y=y_0\in\bar H_2(\bar W:\RR)$. Using the commutative diagram
\eqref{poincare} and by the fact $D(c_1(W))=[V]$, the element
$D^{-1}(y)\in H^2(\bar W,M:\RR)$ satisfies that $j^*(D^{-1}(y))
=c_1(W)$.

 On the other hand, recall from the argument of \cite{b-e2}
that the integral $\langle {\bar c}_1^2(\bar W),[\bar W]\rangle$
does not depend on the choice of lift $\bar c_1(\bar W)$ to $c_1(W)$
, so we can choose $\bar c_1(\bar W)=D^{-1}(y)\in H^2(\bar W,M:\RR)$
(\cf $\S$ \ref{integrand}). By definition, $y\cdot y= \langle {\bar
c}_1(\bar W)^2,[\bar W]\rangle$ which is an integer by Proposition
\ref{integervalue}. This contradiction proves {\bf Theorem}.

\begin{remark}
Neumann and Reid have shown that if an infranil 3-manifold arises as
a cusp cross-section of a 1-cusped complex hyperbolic 2-manifold,
then the rational Euler number must be 1/3-integral. There are
infranilmanifolds which do not satisfy this condition.
\end{remark}



\begin{thebibliography}{99}

\bibitem{b-he}
O. Biquard and M. Herzlich, \lq\lq A Burns-Epstein invariant for
ACHE $4$-manifolds,\rq\rq \ {\sl Duke Math. Jour.}

\bibitem{b-e1}
D. Burns and C. L. Epstein, \lq\lq A global invariant for
three-dimensional $CR$-structure,\rq\rq \ {\sl Invent. Math.},
(1988), 333-348.

\bibitem{b-e2}
D. Burns and C. L. Epstein, \lq\lq Characteristic numbers of bounded
domains,\rq\rq \ {\sl Acta Math.}, (1990), 29-71.


\bibitem{dekim}
K. Dekimpe, \lq Almost-Bieberbach Goups: Affine and Polynomial
Structures,\rq {\sl  Springer-Verlag, Lecture Notes in Math. 1639,
1996.}

\bibitem{e-p}
E. Falbel and J. R. Parker, \lq\lq The geometry of the
Eisenstein-Picard modular group,\rq\rq Preprint.

\bibitem{go} W. Goldman, \lq Complex Hypwrbolic Geometry,\rq Oxford
Univ. Press (1999).

\bibitem{g-h}
P. Griffths and J. Harris, \lq Principles of Algebraic Geometry
\rq \ {\sl John Wiley \& Sons, Inc.} 1978.

\bibitem{g-p}
N. Gusevskii and J. R. Parker, \lq\lq Representations of free
Fuchsian groups in complex hyperbolic space,\rq\rq \ {\sl Topology
39}, (2000), 33-60.

\bibitem{h-r} G.C. Hamrick and D.C. Royster, \lq\lq Flat Riemannian
manifolds are boundaries,\rq\rq \ {\sl Invent. Math.} vol.~66 (2),
405-413, 1982.

\bibitem{h-m}
D. Husemoller and J. Milnor, \lq Symmetric Bilinear Forms,\rq \
 {\sl Springer, Berlin}  1973.


\bibitem{ka3} Y. Kamishima,
\lq\lq Cusp cross-sections of hyperbolic orbifolds by Heisenberg
nilmanifolds I,\rq\rq \ {\sl to appear in Geom. Dedicata}.

\bibitem{k-l-r}
Y. Kamishima, K.B. Lee and F. Raymond, \lq\lq  The Seifert
construction and its application to infranilmanifolds," \ {\sl
Quart. J. Math. Oxford (2)} vol. 34 (1983), 433-452.


\bibitem{ka-tsu}
Y. Kamishima and T. Tsuboi, \lq\lq $CR$-structures on Seifert
manifolds,\rq\rq \ {\sl Invent. Math.}, vol.~104 (1991), pp.
149--163.

\bibitem{l-r}
D.D. Long and A.W. Reid, \lq\lq All flat manifolds are cusps of
hyperbolic orbifolds", \ {\sl Algebraic and Geometric Topology}.
vol. 2 (2002) 285-296.

\bibitem{l-r1}
D.D. Long and A.W. Reid, \lq\lq On the geometric boundaries of
hyperbolic $4$-manifolds", \ {\sl Geometry and Topology}. vol. 4
(2000) 171-178.

\bibitem{ben}
D. B. Mcreynolds, \lq\lq Peripheral separability and cusps of
arithmetic hyperbolic orbifolds",
 \ {\sl Algebr. and Geom. Topol. 2004 (4), 721-755}.


\bibitem{to}
D. Toledo, \lq\lq Representations of surface groups in complex
hyperbolic space,\rq\rq \ {\sl Jour. of Diff. Geom. 29}, (1989),
125-133.





\end{thebibliography}
\end{document}